# ELFPIE: an error-laxity Fourier ptychographic iterative engine

Shuhe Zhang,[1,2] Tos T. J. M. Berendschot,[1] and Jinhua Zhou[2,3]

1. *University Eye Clinic Maastricht, Maastricht University Medical Center + , P.O. Box 5800, Maastricht, AZ 6202, the Netherlands*
2. *School of Biomedical Engineering, Anhui Medical University, Hefei 230032, China.*
3. *Anhui Provincial Institute of Translational Medicine, Anhui Medical University, Hefei 230032, China.*

shuhe.zhang@maastrichtuniversity.nl,
zhoujinhua@ahmu.edu.cn

**Abstract:** We present a simple but efficient and robust reconstruction algorithm for Fourier ptychographic microscopy, termed error-laxity Fourier ptychographic iterative engine (Elfpie), that **is simultaneously robust to (1) noise signal** (including Gaussian, Poisson, and salt & pepper noises), **(2) problematic background illumination problem, (3) vignetting effects,** and **(4) misaligning of LED positions,** without the need of calibrating or recovering of these system errors. In Elfpie, we embedded the inverse problem of FPM under the framework of feature extraction/recovering and proposed a new data fidelity cost function regularized by the global second-order total-variation regularization (Hessian regularization). The closed-form complex gradient for the cost function is derived and is back-propagated using the AdaBelief optimizer with an adaptive learning rate to update the entire Fourier spectrum of the sample and system pupil function. The Elfpie is tested on both simulation data and experimental data and is compared against the state-of-the-art (SOTA) algorithm. Results show the superiority of the Elfpie among other SOTA methods, in both reconstruction quality under different degeneration issues, and implementation efficiency. In general, compared against SOTA methods, the Elfpie is robust to Gaussian noise with 100 times larger noise strength, salt & pepper noise with 1000 times larger noise strength, and Poisson noise with 10 times noise strength. The Elfpie is able to reconstruct high-fidelity sample field under large LED position misalignments up to 2 mm. It can also bypass the vignetting effect in which all SOTA methods fail to reconstruct the sample pattern. With parallel computation, the Elfpie is able to be $K$ times faster than traditional FPM, where $K$ is the number of used LEDs.

## 1. Introduction

Fourier ptychographic microscopy (FPM) is a computational imaging technique that combines concepts of coded-illumination and synthetic aperture [1], offering an effective way to improve the resolution of an imaging system and retrieve the lost phase information of the sample. As such, the FPM attracts a lot of interest in applications including reflective FPM [2, 3], near-field FPM [4], and tomographic FPM (3DFPM) [5-8].

The FPM routines employ either angularly varying illumination (provided by an LED array or a digital mirror device) [1, 9] or aperture scanning (produced by a spatial light modulator) [10, 11] to obtain a series of low-resolution images. Each low-resolution image corresponds to a specific region in the sample's Fourier spectrum with sufficient overlap to provide data redundancy. The low-resolution images are given to the Fourier ptychographic iterative engine (FPIE) for the super-resolution reconstruction which can exceed the diffraction limit of the objective lens. Apparently, the reconstruction quality is related to the raw data quality and the precision of system parameters including LED position, and system aberrations [12, 13].

Since the FPM raw data is sequentially collected, every single low-resolution image is corrupted by unavoidable random noise signal and may also suffer from uneven illumination fluctuations, both can severely degenerate the FPM reconstruction quality. During the reconstruction process, the misaligning of illumination position can lead to problematic spectrum recovery, causing a stripe effect on the reconstructed image. In previous reports, noise reduction, intensity fluctuation correction, and illumination position misaligning are treated separately and have their special designed treatments and assumptions. A direct combination of these algorithms will result in an extremely redundant, unstable, and inefficient FPIE algorithm, which may be of no practical value.

In this research, we re-design the cost function of FPIE including the data fidelity term and penalty term and propose an error-laxity Fourier ptychographic iterative engine, termed Elfpie that is simultaneously robust to noise signal, illumination fluctuation, vignetting effects, and misaligning of LED positions, without calibrating or recovering of these parameters. The subsequent manuscript is organized as follows. Section 2 briefly reviews related works. Section 3 introduces our new cost function under the framework of feature extraction/recovering and presents the core idea of Elfpie and the updating routines. Section 4 presents the simulation study for Elfpie under different degrading conditions. Section 5 presents the experimental results. Section 6 and Section 7 are the discussion and concluding remarks, respectively.

## 2. Related works

### 2.1 Intensity correction

A LEDs array is the widely used in FPM experiments. The illumination intensities of different LED elements behave differently over time due to heating, and voltage changes. In [14], the illumination fluctuation is correct by minimizing the difference in the intensity between calculated low-res images and corresponding measured low-res images. A certain ratio is calculated and assigned to the measured low-res images to compensate for the illumination fluctuation. However, this method cannot correct the background intensity fluctuation caused by interference of inhomogeneous media distribution or shadows cast by other defocusing layers.

To tackle the background intensity fluctuation problems, an adaptive background correction [15] is proposed by assigning additional background terms to each raw image which can be recovered during the FPIE iteration. Another intensity correction involves creating a binary mask according to a threshold to remove the problematic area from current iterations [16, 17]. Although the binary mask causes loss of image information, the reconstruction quality will not decrease a lot as the raw data contains data redundancy.

Furthermore, the FPM data may also suffer from vignetting effect [18, 19], where some images contain both partially bright- and partially dark-field information. The vignetting effect results from the non-linear effect of a practical image system. In traditional FPM, there is no way to overcome the vignetting effect and people often segment the image into bright and dark parts or skip the images containing vignetting in order to obtain good FPM recovery results.

### 2.2 LED position correction

The misaligning of LED positions will lead to stripe-like artifacts [12] as the recovered spectrum components are placed in the wrong spectral positions, which further can significantly decrease the FPM reconstruction. The LED misaligning can always happen to FPM experiments with different misaligning degrees. Correction of LED position can promote a high-quality FPM reconstruction.

The LED position correction can be categorized into individual correction and global correction strategy, based on the modeling of LED position shifting. In the individual strategy, the position of each LED dot has its own random shifting and is separately corrected by searching for the possible positions where the cost function is minimal. Using the simulated annealing algorithm is feasible but is very time-consuming since a 2D-plane search is needed for every single LED dot [12, 20]. The search range is also dependent on the distance of misaligning. Later gradient-descent methods [21, 22] were proposed by incorporating the LED position into the system parameters to be optimized. The LED correction can be also achieved as a preprocessing step before starting an FPM reconstruction using the self-calibration method [23].

In the global correction strategy, the relative position between each LED dot is assumed to be unchanged as they embedded in the LED panel. While the entire LED panel can have a certain transverse shifting and may also rotate in some degrees [13, 24]. By using a simulated annealing algorithm, both the offset and the rotation angle can be found during the FPM reconstruction. Compared to the individual strategy, the global one has less system parameters to be optimized but cannot be applied to more general cases.

In [25] a multi-look method is introduced by digitally changing the system parameters and calculating the reconstruction using different parameters to obtain a series of reconstruction results. The final output is the average of them, yielding a miscalibration-tolerant FPM reconstruction.

### 2.3 Noise suppression

As many dark-field images are being captured during the FPM data acquisition, the noise signals will largely decrease the reconstruction quality as the signal-to-noise ratio for the dark-field images is rather lower than those of bright-field images.

To tackle the impact of the noise signal, many researches have been proposed, including the Wirtinger flow (WF) optimization [26], which put the FPM into the Wirtinger flow framework and is robust to Gaussian noise; the adaptive step-size (AP) strategy [27], in which the step size for the gradient descent is changed according to the estimation of residual; and threshold-based methods [17, 28, 29] methods, where a threshold is used to erase pixels that are regarded as noises in the dark-field images. It is shown that the WF method works well for Gaussian noise, while the AS method performs better than the AP method for Poisson noise. The threshold method works well for mild noise and will cause a loss of reconstruction quality when the noise becomes noticeable. Further, the FPM is combined with total variation (TV) regularization to suppress the noise signal. In this study, the TV regularization is treated as an additional module that is applied to the amplitude and phase components during each iteration [30, 31].

*2.4 Our contribution*

Despite the simplicity of the FPM system, a robust FPM reconstruction relies on several conditions including stable illumination/background intensity, the high signal-to-noise ratio of raw low-res images, and precise determination of illumination direction. Each of them has its solution which is independent to the others. According to our investigation, there is no FPM algorithm that can simultaneously solve all these three problems. In this research, we redesign the FPIE and proposed the Elfpie for experimental robust FPM reconstruction. The Elfpie has the following contributions:

(1) We proposed a new cost function for FPM reconstruction. Instead of using the amplitude of low-res image for FPM reconstruction, we use the edge information, which is the first-order spatial gradient of the image for FPM reconstruction and applied the $L_1$-norm to the data fidelity term. By doing this, **the Elfpie can bypass the impact of illumination/background fluctuations, and vignetting effects, and is robust to misaligning of LED position and potential salt & pepper (SNP) noise.**

(2) We introduced the second-order TV regularization (Hessian regularization) as our penalty function to regularize the second-order spatial gradient of both amplitude and phase components of the recovered high-res image. By doing this, **the Gaussian, Poisson, and SNP noises can be efficiently suppressed.**

(3) We calculated the closed-form complex gradient of our cost function; the gradient is back-propagated using the AdaBelif optimizer with an adaptive learning rate to update the sample's Fourier spectrum. With that, **the algorithm can converge fast and has a high probability to avoid the local minimum points.**

(4) Finally, our Elfpie can be parallel as the Fourier spectrum region corresponding to each LED can be updated simultaneously, which, theoretically, can be ***K* times (*K* is the number of LEDs) faster than the traditional FPM reconstruction routine.**

## 3. Error-laxity Fourier ptychographic iterative engine

*3.1 FPM forward model*

We consider the framework of multiplexed FPM [32, 33] as a universal framework for FPM raw data acquisition routine, in which the LEDs are assigned into *N* groups, and each group contains *M* LEDs. *M* is a factor of *N* so that every LED is only lighted once. The individual raw low-res image is captured when the sample is illuminated by *M* LEDs simultaneously, therefore, a total of *N* raw images is obtained.

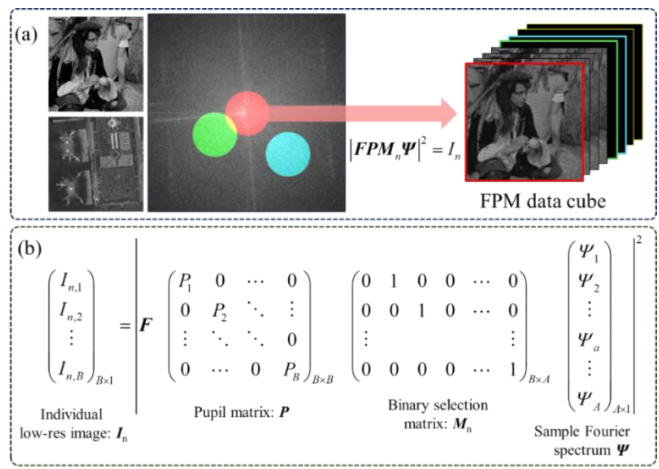

Fig. 1. Forward model of FPM. (a) sketch of image formation process in FPM. (b) matrix representation.

Let $\boldsymbol{O} = \mathcal{A}\exp(i\boldsymbol{\Phi})$, $\boldsymbol{O} \in \mathbb{C}^A$ be the complex wavefront of the sample (a total of *A* pixels) and $\boldsymbol{\Psi} = \boldsymbol{F}_A^{-1}\boldsymbol{O}$ be the Fourier spectrum of vectorized $\boldsymbol{O}$, where $\boldsymbol{F}_A$, $\mathbb{C}^A \to \mathbb{C}^A$ denotes the discrete Fourier transform matrix. The forward model of FPM in an ideal condition can be described as [32]

$$\hat{\boldsymbol{I}}_n(\boldsymbol{r}) = \sum_{m=1}^{M} \left| \boldsymbol{F}_B \boldsymbol{P} \boldsymbol{M}_{n,m} \boldsymbol{\Psi} \right|^2, \quad \hat{\boldsymbol{I}}_n \in \mathbb{R}^B, \quad \mathbb{C}^A \to \mathbb{R}^B, \tag{1}$$

where $\boldsymbol{r}$ is the spatial coordinate vector. $\boldsymbol{P}$, $\mathbb{C}^B \to \mathbb{C}^B$ is a diagonal matrix denoting the pupil function of the objective lens. $\boldsymbol{M}_{n,m}$, $\mathbb{C}^A \to \mathbb{C}^B$ denotes the illumination matrix which is a real-value binary mask that extracts a certain part of the regions (a total of *B* pixels, $B < A$) of $\boldsymbol{\Psi}$ that corresponds to the *m*-th LED dot in *n*-th group. $\boldsymbol{M}_{n,m}^T$, $\mathbb{C}^B \to \mathbb{C}^A$ denotes an inverse operation that put the extracted vector into the origin place. The data acquisition routine reduces to traditional sequential FPM when $M = 1$.

As mentioned in section 2, the image can be degenerated by uneven illumination fluctuations, background intensity, and noises. The forward model in Eq. (1) can be rewritten for a more general case which is

$$\boldsymbol{I}_n(\boldsymbol{r}) = \sum_{m=1}^{M} \Lambda_m(\boldsymbol{r}) \circ \left| \boldsymbol{F}_B \boldsymbol{P} \boldsymbol{M}_{n,m} \boldsymbol{\Psi} \right|^2 + \boldsymbol{\Omega}_n + \boldsymbol{\xi}_n, \quad \mathbb{C}^A \to \mathbb{R}^B, \tag{2}$$

$\Lambda_m$ and $\boldsymbol{\xi}_n$ are degeneration agents, where $\Lambda_m$ denotes the illumination fluctuations, $\boldsymbol{\Omega}_n$ is the additive background intensity and $\boldsymbol{\xi}_n$ is the noise signal.

The FPM recovery is a maximum a posteriori estimate (MAP) problem [34-36], in which we wish to reconstruct $\boldsymbol{\Psi}$ under the series of observations of $\boldsymbol{I}_n$, $n \in [1, N]$, and can be solved using a complex least-square method with gradient descent. However, the reconstruction can be severely degraded as the quality of $\boldsymbol{I}_n$ is rather lower than that of $\hat{\boldsymbol{I}}_n$ due to the presence of degeneration agents, for example noise signals. Since they may change rapidly during the data acquisition, manual calibration can be time-consuming and unstable. We hope to suppress or bypass their impact (at least for $\Lambda_m$ and $\boldsymbol{\Omega}_n$) during the reconstruction process. To do so, we introduce the idea of feature extraction/recovering.

### 3.2 Feature extraction/recovering

Imaging that if both $\hat{\boldsymbol{I}}_n$ and $\boldsymbol{I}_n$ contain their comment features which are unique and won't present in $\Lambda_m$ and $\boldsymbol{\Omega}_n$, we are able to use these features for FPM reconstruction and bypass the impact of $\Lambda_m$. Let $\boldsymbol{\Theta}$ being an invertible operator that can extract such unique features. Applying $\boldsymbol{\Theta}$ to $\hat{\boldsymbol{I}}_n$ and $\boldsymbol{I}_n$ yields $\boldsymbol{\Theta}\hat{\boldsymbol{I}}_n \simeq \boldsymbol{\Theta}\boldsymbol{I}_n$. By solving the following optimal function

$$(\boldsymbol{\Psi}, \boldsymbol{P}) = \arg\min \left\langle \boldsymbol{\Theta} g(\boldsymbol{I}_n), \; \boldsymbol{\Theta} g\left( \sum_{m=1}^{M} \left| \boldsymbol{F}_B \boldsymbol{P} \boldsymbol{M}_{n,m} \boldsymbol{\Psi} \right|^2 \right) \right\rangle, \tag{3}$$

we are able to obtain a reconstructed $\boldsymbol{\Psi}$ with is less impacted by $\Lambda_m$ and $\boldsymbol{\Omega}_n$. Here $\langle X, Y \rangle$ denotes an arbitrary function that can measure the "distance" between variables *X* and *Y*. (for example the Euclidean distance, also known as $L_2$-norm). *g* is a pixel-wise scaling function that re-maps the intensity measurement to a certain distribution. For example, $g(x) = x$ denoting using the intensity measurement, while $g(x) = \sqrt{x}$ denoting using the amplitude measurement.

Accordingly, the quality of $\boldsymbol{\Psi}$ now depends on the selection of $\boldsymbol{\Theta}$, *g*, and the distance function that $\boldsymbol{\Theta}$ should efficiently separate the image features from the degeneration agents, while g and distance function should maintain the data fidelity for successful reconstruction results. The $\boldsymbol{\Theta}$, *g*, and the distance function can be either learned by using neural networks [37], or manually chosen according to the statistical features of the captured images [38]. Considering that both $\Lambda_m$ and $\boldsymbol{\Omega}_n$ do not contain edge information as they originate from illumination problems (spatially slow varying), we choose the spatial gradient operator to extract the edge feature of the image for the reconstruction of $\boldsymbol{\Psi}$. Here $\boldsymbol{\Theta} = \nabla$ is the spatial gradient operator. The cost function for the amplitude-based data fidelity, where $g(x) = \sqrt{x}$, is given by

$$\mathcal{L}_{\text{Fidelty-Amplitude}}(\boldsymbol{\Psi},\boldsymbol{P}) = \sum_{n=1}^{N}\left\|\nabla\sqrt{\boldsymbol{I}_n} - \nabla\sqrt{\sum_{m=1}^{M}|\boldsymbol{F}_B\boldsymbol{P}\boldsymbol{M}_{n,m}\boldsymbol{\Psi}|^2}\right\|_1, \quad \mathbb{C}^A \to \mathbb{R}. \tag{4}$$

While the intensity-based data fidelity when $g(x)=x$ is

$$\mathcal{L}_{\text{Fidelty-Intensity}}(\boldsymbol{\Psi},\boldsymbol{P}) = \sum_{n=1}^{N}\left\|\nabla\boldsymbol{I}_n - \nabla\sum_{m=1}^{M}|\boldsymbol{F}_B\boldsymbol{P}\boldsymbol{M}_{n,m}\boldsymbol{\Psi}|^2\right\|_1, \quad \mathbb{C}^A \to \mathbb{R}. \tag{5}$$

Instead of using L$_2$-norm, here we use the L$_1$-norm to promote the sparsity of the gradient of the illumination component in the recovered $\boldsymbol{O}$. As such, the background in the recovered image can be more uniformly distributed. Moreover, the introduction of gradient operation allows our Elfpie to bypass the impact of the mismatching of low-frequency component caused by LED misaligning.

To further suppress the noise signal, we chose the second-order total variation, the Hessian variation, as the penalty function to suppress the second-order spatial gradient of both amplitudes $\mathcal{A}$:

$$\mathcal{L}_{\mathcal{A}\text{-Hessian}}(\boldsymbol{\Psi}) = \left\|\boldsymbol{H}|\boldsymbol{O}|\right\|_1 = \left\|\boldsymbol{H}\sqrt{\text{Re}^2\boldsymbol{F}_A\boldsymbol{\Psi}+\text{Im}^2\boldsymbol{F}_A\boldsymbol{\Psi}}\right\|_1, \quad \mathbb{C}^A \to \mathbb{R}, \tag{6}$$

and phase $\boldsymbol{\Phi}$:

$$\mathcal{L}_{\boldsymbol{\Phi}\text{-Hessian}}(\boldsymbol{\Psi}) = \left\|\boldsymbol{H}\angle\boldsymbol{O}\right\|_1 = \left\|\boldsymbol{H}\arctan\frac{\text{Im}\,\boldsymbol{F}_A\boldsymbol{\Psi}}{\text{Re}\,\boldsymbol{F}_A\boldsymbol{\Psi}}\right\|_1, \quad \mathbb{C}^A \to \mathbb{R}, \tag{7}$$

Here $\boldsymbol{H}=\nabla\nabla$ is the spatial Hessian operator calculating the second-order gradient in $x$, $y$, and $x$-$y$ directions. Instead of using the total-variation regularization which can cause staircase-like effects, the Hessian regularization obtains more smoothly varying pixels.

The global cost function for the Elfpie is the linear combination of Eq. (4) [or (5)] to Eq. (7) which is given by

$$\mathcal{L}(\boldsymbol{\Psi},\boldsymbol{P}) = \mathcal{L}_{\text{Fidelity}}(\boldsymbol{\Psi},\boldsymbol{P}) + \alpha\mathcal{L}_{\mathcal{A}\text{-Hessian}}(\boldsymbol{\Psi}) + \beta\mathcal{L}_{\boldsymbol{\Phi}\text{-Hessian}}(\boldsymbol{\Psi}), \tag{8}$$

where $\alpha$ and $\beta$ are the penalty parameters. $\boldsymbol{\Psi}$ can be obtained by minimizing the cost function. To do this, we need to calculate the gradients of Eq. (8) and back-propagate it to update $\boldsymbol{\Psi}$. Unlike traditional FPIE where the Fourier spectrum is updated part by part through a certain sequence (spiral, or row-by-row), our Elfpie updates the entire Fourier spectrum in each iteration.

*3.3 Learning for ψ*

Eq. (8) involves three L$_1$-norm terms that are non-differentiable at the origin. While in this case, the cost function is zero and there is no need to update the parameters. By assuming the gradient of L$_1$-norm at origin to be zero, and using the $\mathbb{CR}$-calculus (known also as Wirtinger calculus) [39], the closed-form gradient of $\mathcal{L}_{\text{Fidelity}}$ with respect to $\overline{\boldsymbol{\Psi}}$ (complex conjugate of $\boldsymbol{\Psi}$) is written as

$$\nabla_{\overline{\boldsymbol{\Psi}}}\mathcal{L}_{\text{Fidelity}} = \sum_{n=1}^{N}\sum_{m=1}^{M}\boldsymbol{M}_{n,m}^{\mathcal{H}}\boldsymbol{P}^{\mathcal{H}}\boldsymbol{F}_B^{\mathcal{H}}\boldsymbol{W}_{n,m}. \tag{9}$$

where

$$\boldsymbol{W}_{n,m} = \text{Diag}\left[\boldsymbol{o}_{n,m}\circ\left(\sum_{m=1}^{M}|\boldsymbol{o}_{n,m}|^2\right)^{-1/2}\right]\nabla^T\omega\left(\sqrt{\sum_{m=1}^{M}|\boldsymbol{o}_{n,m}|^2}-\sqrt{\boldsymbol{I}_n};\nabla\right), \tag{10}$$

calculated from Eq. (4) (Amplitude-based), and

$$\boldsymbol{W}_{n,m} = \text{Diag}(\boldsymbol{o}_{n,m})\nabla^T\omega\left(\sqrt{\sum_{m=1}^{M}|\boldsymbol{o}_{n,m}|^2}-\sqrt{\boldsymbol{I}_n};\nabla\right), \tag{11}$$

calculated from Eq. (5) (Intensity-based). $\boldsymbol{o}_{n,m}=\boldsymbol{F}_B\boldsymbol{P}\boldsymbol{M}_{n,m}\boldsymbol{\Psi}$, and the superscript $\mathcal{H}$ denotes the Hermitian transpose and $T$ denotes the transpose. $\circ$ is the element-wise product. $\omega(X;\nabla)$ is the normalization function that calculates a certain result between operator $\nabla$ and variable $X$ in an isotropic or an anisotropic manner since $\nabla=(\nabla_x,\nabla_y)^T$ has two components. $\omega(X;\nabla)$ is defined as:

$$\omega(X;\nabla) = \begin{cases} \dfrac{1}{\sqrt{(\nabla_x X)^2 + (\nabla_y X)^2}} (\nabla_x X, \nabla_y X)^T, & \text{For isotropic} \\ (\text{sign}(\nabla_x X), \text{sign}(\nabla_y X))^T, & \text{For an-isotropic} \end{cases}. \quad (12)$$

The gradient of Eq. (5) and Eq. (6) with respect to $\Psi$ is given by

$$\nabla_{\bar{\Psi}} \mathcal{L}_{\mathcal{A}\text{-Hessian}} = \frac{1}{2} F_A^{\mathcal{H}} \text{Diag}\left(\frac{O}{|O|}\right) H^T \omega(|O|; H), \quad (13)$$

and

$$\nabla_{\bar{\Psi}} \mathcal{L}_{\Phi\text{-Hessian}} = \frac{i}{2} F_A^{\mathcal{H}} \text{Diag}\left(\frac{O}{|O|^2 + \eta}\right) H^T \omega(\angle O; H), \quad (14)$$

respectively, where $\eta$ is a small value to avoid dividing by zeros condition. The $\omega(|O|; H)$ and $\omega(\angle O; H)$ can be also calculated in an isotropic or an anisotropic manner.

The gradient of the global cost function is then given by

$$\nabla_{\bar{\Psi}} \mathcal{L} = \nabla_{\bar{\Psi}} \mathcal{L}_{\text{Fidelity}} + \alpha \cdot \nabla_{\bar{\Psi}} \mathcal{L}_{\mathcal{A}\text{-Hessian}} + \beta \cdot \nabla_{\bar{\Psi}} \mathcal{L}_{\Phi\text{-Hessian}}, \quad (15)$$

which is backpropagated to update the parameter $\Psi$.

The parameters $\alpha$ and $\beta$ are automatically determined based on

$$\alpha = \beta = \frac{1}{5} \sqrt{\frac{\pi}{2}} \frac{1}{N} \frac{1}{WH} \sum_{n=1}^{N} \sum_{x=1}^{W} \sum_{y=1}^{H} |I_n(x,y) \otimes \mathcal{L}|, \quad (16)$$

for intensity-based Elfpie, and

$$\alpha = \beta = \frac{1}{5} \sqrt{\frac{\pi}{2}} \frac{1}{N} \frac{1}{WH} \sum_{n=1}^{N} \sum_{x=1}^{W} \sum_{y=1}^{H} \left|\sqrt{I_n(x,y)} \otimes \mathcal{L}\right|, \quad (17)$$

for amplitude-based Elfpie. $\mathcal{L} = [-1, 2, -1; -2, 4, -2; -1, 2, -1]$ is the Laplacian operator.

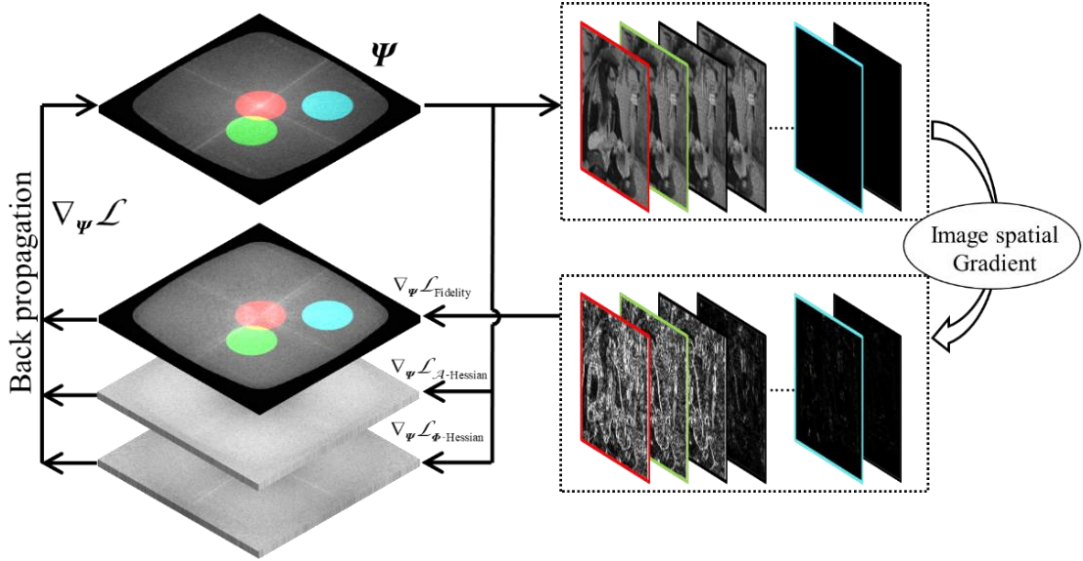

Fig. 2. Backpropagation of gradient in Elfpie for $\Psi$ updating.

Fig. 2 shows the backpropagation of gradient of learning for $\Psi$. With $\Psi$ from last step, we obtain the calculated low-resolution images and obtain their spatial gradients. Then we use Eq. (10) or Eq. (11) to obtain the gradient component, $\nabla_{\Psi} \mathcal{L}_{\text{Fidelity}}$, from data fidelity term. The $\nabla_{\Psi} \mathcal{L}_{\mathcal{A}\text{-Hessian}}$ and $\nabla_{\Psi} \mathcal{L}_{\Phi\text{-Hessian}}$ from amplitude and phase Hessian regularization are directly calculated from the $\Psi$ using Eq. (13) and Eq. (14). Then the entire Fourier spectrum is updated.

### 3.4 Learning for P

The pupil function $P$ can be also learned from the data fidelity term. The gradient of Eq. (4) and Eq. (5) with respect to $P$ is

$$\nabla_{\bar{P}} \mathcal{L}_{\text{Fidelity}}(\Psi, P) = \sum_{n=1}^{N} \sum_{m=1}^{M} \text{Diag}(M_{n,m}\Psi)^{\mathcal{H}} F^{\mathcal{H}} W_{n,m}, \qquad (18)$$

Where $W_{n,m}$ is calculated from Eq. (10) for amplitude-based fidelity, and Eq. (11) for intensity-based fidelity, respectively. Parameter $P$ can be further learned from the gradient.

It is worth noting that the global cost function is a real-value function of multi-complex variables. The cost function is non-convex and using the convex optimization method may lead the solution convergent to a local minimum. To accelerate the gradient descent and give the possibility to escape from the local minimum, here we propose a modified AdaBelif optimizer [40] with an adaptive learning rate [41] to update the model parameters.

$$\begin{cases}
g = \nabla_{\Psi} \mathcal{L}(\Psi^{t-1}) & \\
\mu^t = \gamma_1 \cdot \mu^{t-1} + (1-\gamma_1) \cdot g & \text{(First moment estimate)} \\
v^t = \gamma_2 \cdot v^{t-1} + (1-\gamma_2) \cdot |\mu^t - g|^2 & \text{(Belief in gradient direction)} \\
\hat{\mu} = \mu^t / (1-\gamma_1^t) & \\
\hat{v} = v^t / (1-\gamma_2^t) & \text{(Bias correction)} \\
\Delta\Psi = \frac{\sqrt{\delta^{t-1}} + \eta}{\sqrt{\hat{v}} + \eta} \circ [\gamma_1 \cdot \hat{\mu} + (1-\gamma_1) \cdot g] & \text{(Increment of } \Psi) \\
\Psi^t = \Psi^{t-1} - \Delta\Psi & \\
\delta^t = \gamma_1 \cdot \delta^{t-1} + (1-\gamma_1) \cdot |\Delta\Psi|^2 & \text{(Adaptive learning rate)}
\end{cases} \qquad (19)$$

At the $t$-th iteration, the $\Psi$ is updated according to Eq. (19), and $\mu^0 = v^0 = 0$, $\delta^0 = 1$, and $r_1 = 0.9$, $r_2 = 0.999$ by default. $\Delta\Psi$ denotes the increment of $\Psi$ in each iteration, and the learning rate $\sqrt{\delta}$ is adaptively determined according to the historical increment of $\Psi$ ($\Delta\Psi$).

Our proposed optimizer allows dynamic adjustment of the learning rate, the model won't stay in saddle point for long and has the ability to escape from the local minima, which convergent fast (within 20 iterations). In general, the Elfpie is summarized in **Algorithm 1**. Accordingly, our Elf-pie updates the entire Fourier spectrum of the sample once in each outer iteration. This feature is different from the traditional FPIE where parts of the spectrum are updated in each inner iteration (iteration for LED dots). Since the calculation of $W_{n,m}$ and $V_{n,m}$ is independent to each other, the Elfpie is suitable for parallel computation, which can be $N*M$ times faster than traditional FPIE.

The physical insight of Eq. (10) is obvious when $M = 1$. Instead of replacing the calculated amplitude $|o_n|$ by the measured amplitude $\sqrt{I_n}$, in the Elfpie, the amplitude is replaced by the difference between the gradient of $|o_n|$ and $\sqrt{I_n}$ with a certain normalization function $\omega(X; \nabla)$. In such a manner, we are able to avoid the impact of low-frequency problematic illumination problems, making the recovery algorithm more robust to illumination problem and vignetting effect than traditional FPIE.

**Algorithm 1: Elf-pie updating algorithm**

**Input**: FPM raw data cube, $I_n$, parameters $\alpha$, $\beta$, and $\gamma_1$, $\gamma_2$. $\boldsymbol{\mu}^0 = \boldsymbol{v}^0 = 0$

*For* $t = 1, 2, \cdots$ iter_max
    *For* $n = 1, 2, \cdots N$
        *For* $m = 1, 2, \cdots M$
            Calculate $W_{n,m}$ using Eq. (10) or Eq. (11)    } can be parallel
        *End For*
    *End For*
    Calculating $\nabla_{\boldsymbol{\Psi}} \mathcal{L}_{\text{Fidelity}}$ and $\nabla_P \mathcal{L}_{\text{Fidelity}}$ using Eq. (9) and Eq. (18)
    $O \leftarrow F_A^{-1} \boldsymbol{\Psi}$
    Calculating $\nabla_{\boldsymbol{\Psi}} \mathcal{L}_{A\text{-Hessian}}$ using Eq. (13)
    Calculating $\nabla_{\boldsymbol{\Psi}} \mathcal{L}_{\Phi\text{-Hessian}}$ using Eq. (14)
    Learning for $\boldsymbol{\Psi}$ using Eq. (19)
    Learning for $P$ using Eq. (19)
*End For*

**Output**: High resolution complex wavefront $O$

Note: Before updating $P$, the resulting $P$ from the last step is filtered by a Gaussian kernel with a size of 3 by 3 pixels, and a width of 1 pixel to avoid unexpected pulse error that normally happens in real experiments.

## 4. Simulation study

### 4.1 Experimental setting

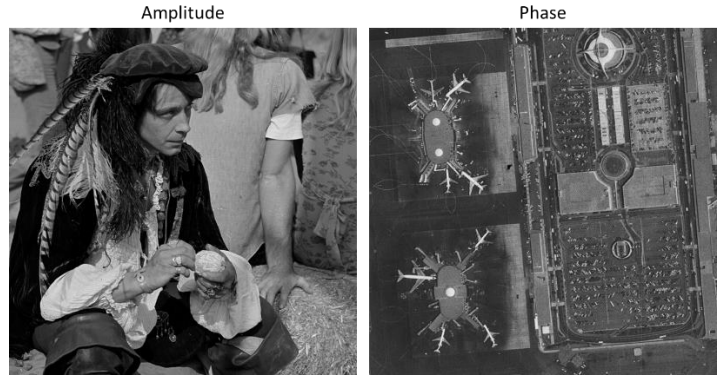

Fig. 3, Amplitude (left) and phase (right) pattern used in the simulation study.

In the simulation section, we use the images shown in Fig 3 for the amplitude component and phase component as the complex wavefront of the sample. The grayscale for both images is normalized between (0.1 and 1) in order to avoid the condition where some pixels have zero intensity. The high-resolution ground truth images are 513 by 513 pixels, and the simulated images captured by a camera are 128 by 128 pixels. The NA of the objective lens is 0.10, and the magnification is ×4. The pixel size of the camera is 3.65 μm. The sample is illuminated by a LED array of 15 by 15 LED dots with a 536 nm wavelength. The array is placed 90 mm away from the sample, and the distance between two adjacent LED dots is 6 mm. In an ideal condition, the very center LED dots is placed at the optics axis. No LED multiplexed is used.

The program is composed in MATLAB 2022 and we use the "*randn*" command to add Gaussian noise, and "*imnoise*" command to add SNP noise to the image. For Gaussian, and SNP noise, the noise strength can be controlled by the variation value $a$. The method of adding Gaussian or SNP noise at different noise levels is straightforward while adding Poisson noise is more complex as it relates to the local intensity level of the input image. In our research, we add Poisson noise to the image by simulating the photon counting process, denoted by PN in Eq. (20). In our research, we refer the Poisson noise level by Lv1 to Lv4 from weak to strong noise.

Furthermore, we use the convolution pattern between a Gaussian noise, $\zeta_n$, (a = 0.001) and a Gaussian blur kernel, $G$, 15 pixels half-waist to simulate the uneven problematic illumination pattern. In general, the forward model in our simulation can be described as

$$\begin{cases} \boldsymbol{I}_n = \left[(1-c)+c\cdot\boldsymbol{G}\otimes\boldsymbol{\zeta}_n\right]\circ\left|\boldsymbol{F}_B\boldsymbol{PM}_n\boldsymbol{\Psi}\right|^2 + a\times\text{randn}(B), & \text{Gaussian noise} \\ \boldsymbol{I}_n = imnoise\left(\left[(1-c)+c\cdot\boldsymbol{G}\otimes\boldsymbol{\zeta}_n\right]\circ\left|\boldsymbol{F}_B\boldsymbol{PM}_n\boldsymbol{\Psi}\right|^2, \text{'salt \& pepper'}, a\right), & \text{SNP noise} \\ \boldsymbol{I}_n = \text{PN}\left(\left[(1-c)+c\cdot\boldsymbol{G}\otimes\boldsymbol{\zeta}_n\right]\circ\left|\boldsymbol{F}_B\boldsymbol{PM}_n\boldsymbol{\Psi}\right|^2\right), & \text{Poisson noise} \end{cases} \quad (20)$$

the parameter $c$ is to control the strength of uneven illumination.

To quantitative measure the image corruption level between idea image and simulated degenerated image, we use [27]

$$NL = \frac{1}{\#DFI}\sum_{\forall all\ DFI}\frac{\|I_{simu}-I_{noised}\|_1}{\|I_{simu}\|_1}\times 100\%. \quad (21)$$

The DFI denotes the dark-field image and $\#DFI$ denotes the number of DFI. The calculation of (21) is based only on DFI since the reconstruction quality of FPM depends more on dark-field images than bright-field images as many dark-field images are collected.

To quantitative measure the reconstruction quality, we use the linear regressed SNR (LSNR) to get rid of additive constants $b$. Between the ground-truth map $\tilde{\varphi}$ and a reconstructed one $\varphi$, the LSNR measure is defined as

$$\text{LSNR}(\varphi,\tilde{\varphi}) = \max_{b\in\mathbb{R}}\ 10\log\left[\frac{\|\tilde{\varphi}\|_2^2}{\|\tilde{\varphi}-(\varphi+b)\|_2^2}\right]. \quad (22)$$

The LSNR is applied to both recovered amplitude and phase. According to Eq. (20), the larger value of LSNR implies better phase reconstruction quality.

We compared our Elfpie against three state-of-the-art methods and they are marked as Adaptive step-size FPIE (AS-FPIE), ADMM-FPIE, and FPIE with momentum (mFPIE). The MATLAB codes for AS-FPIE and can be found in https://www.scilaboratory.com/code.html, while code for mFPIE method and https://github.com/SmartImagingLabUConn/Fourier-Ptychography, respectively. The model parameters are given by default in the codes.

*4.2 Simulation results*

**Table 1. Noise type and corresponding noise levels used in the simulation study**

| Noise type | | Noise Level | | | |
|---|---|---|---|---|---|
| | | a = 0.00001 | a = 0.0001 | a = 0.001 | a = 0.01 |
| Gaussian noise | C=0.25 | 20.49 % | 99.18 % | 800.10 % | 7669.71 % |
| | C=0.50 | 31.02 % | 102.59 % | 796.48 % | 7662.03 % |
| | C=0.75 | 41.98 % | 107.06 % | 793.00 % | 7656.55 % |
| | | a = 0.001 | a = 0.01 | a = 0.1 | a = 0.2 |
| SNP noise | C=0.25 | 989.06 % | 9583.58 % | 47757.62 % | 95453.10 % |
| | C=0.50 | 971.87 % | 9648.64 % | 47799.04 % | 95697.25 % |
| | C=0.75 | 972.56 % | 9543.86 % | 47771.76 % | 95571.99 % |
| | | Lv 1 | Lv 2 | Lv 3 | Lv 4 |
| Poisson noise | C=0.25 | 36.37 % | 66.13 % | 81.07 % | 94.23 % |
| | C=0.50 | 41.78 % | 67.58 % | 80.57 % | 92.12 % |
| | C=0.75 | 49.39 % | 70.74 % | 81.67 % | 91.36 % |

The noise level and background intensity fluctuation used in this research are listed in Tab. 1. We test the methods with different noise types, noise levels, intensity fluctuation, and LED shifting. Each group is repeated 10 times, and the average LSNR is calculated. The mean value of LSNR for the amplitude component and phase component is the final LSNR score for each routine.

Table 2. Reconstruction quality of different methods for Gaussian noise.

| LED position shifting $d$ (mm) | Uneven illumination level $c$ | Noise level $a$ | Average LSNR for different methods, **Gaussian noise** | | | | |
|---|---|---|---|---|---|---|---|
| | | | mFPIE | ADMM-FPIE | AS-FPIE | Elfpie_A | Elfpie_I |
| 0 | 0.25 | $10^{-5}$ | 55.52 | 56.58 | 58.21 | 53.85 | 50.23 |
| | | $10^{-4}$ | 45.43 | 51.97 | 48.92 | 45.02 | 47.50 |
| | | $10^{-3}$ | 28.09 | 44.20 | 30.08 | 46.22 | 47.09 |
| | | $10^{-2}$ | 8.82 | 39.72 | 10.09 | 42.11 | 41.28 |
| | 0.50 | $10^{-5}$ | 44.50 | 50.43 | 50.12 | 49.42 | 47.77 |
| | | $10^{-4}$ | 40.64 | 48.64 | 44.84 | 48.28 | 47.01 |
| | | $10^{-3}$ | 26.16 | 41.96 | 28.39 | 42.48 | 42.18 |
| | | $10^{-2}$ | 7.60 | 38.34 | 8.98 | 39.82 | 39.58 |
| | 0.75 | $10^{-5}$ | 37.29 | 43.03 | 43.17 | 46.81 | 44.20 |
| | | $10^{-4}$ | 33.04 | 41.15 | 38.77 | 39.66 | 43.46 |
| | | $10^{-3}$ | 23.80 | 38.12 | 26.51 | 35.32 | 39.24 |
| | | $10^{-2}$ | 6.27 | 35.49 | 7.75 | 32.88 | 36.10 |
| 0.5 | 0.25 | $10^{-5}$ | 28.74 | 36.64 | 35.89 | 41.87 | 43.47 |
| | | $10^{-4}$ | 15.01 | 26.22 | 20.72 | 38.25 | 41.48 |
| | | $10^{-3}$ | 16.29 | 26.41 | 18.24 | 34.58 | 38.81 |
| | | $10^{-2}$ | 1.56 | 20.18 | 2.57 | 32.72 | 35.40 |
| | 0.50 | $10^{-5}$ | 17.99 | 25.41 | 24.93 | 40.35 | 42.01 |
| | | $10^{-4}$ | 18.81 | 29.13 | 23.96 | 36.90 | 40.59 |
| | | $10^{-3}$ | 11.91 | 21.96 | 13.91 | 34.55 | 37.94 |
| | | $10^{-2}$ | 0.75 | 21.44 | 1.81 | 31.47 | 34.96 |
| | 0.75 | $10^{-5}$ | 19.37 | 28.32 | 26.10 | 39.10 | 39.78 |
| | | $10^{-4}$ | 18.59 | 24.90 | 23.91 | 35.76 | 38.29 |
| | | $10^{-3}$ | 15.71 | 27.86 | 18.46 | 33.42 | 36.64 |
| | | $10^{-2}$ | 2.05 | 24.83 | 3.28 | 31.39 | 34.38 |
| 2 | 0.25 | $10^{-5}$ | 4.19 | 13.14 | 10.89 | 31.51 | 33.40 |
| | | $10^{-4}$ | 0.64 | 8.25 | 3.98 | 31.22 | 32.85 |
| | | $10^{-3}$ | -2.92 | 13.22 | 3.49 | 30.26 | 32.70 |
| | | $10^{-2}$ | -2.86 | 9.66 | -1.39 | 29.50 | 30.45 |
| | 0.50 | $10^{-5}$ | 0.87 | 10.07 | 8.09 | 32.06 | 33.65 |
| | | $10^{-4}$ | -2.05 | 12.12 | 4.73 | 31.45 | 34.48 |
| | | $10^{-3}$ | 0.95 | 13.00 | 3.77 | 30.33 | 32.13 |
| | | $10^{-2}$ | -4.11 | 13.02 | -3.03 | 29.83 | 31.27 |
| | 0.75 | $10^{-5}$ | 1.38 | 10.50 | 6.30 | 31.30 | 33.50 |
| | | $10^{-4}$ | 0.73 | 13.85 | 9.87 | 30.25 | 32.40 |
| | | $10^{-3}$ | -2.38 | 11.15 | 3.16 | 30.36 | 31.65 |
| | | $10^{-2}$ | -5.34 | 9.64 | -3.22 | 29.80 | 29.95 |

Results for Gaussian noise are listed in Tab. 2. When there is no LED shifting, and the noise level and intensity fluctuation are low, all five methods obtain high LSNR values as they are larger than 50. Within, the AS-FPIE method gains the best reconstruction quality. When the intensity fluctuation is 0.75, our proposed Elfpie obtains better LSNR than the other four methods. When the noise level becomes large, the LSNRs for both mFPIE and AS-FPIE decrease significantly, which fail to obtain good reconstruction results.

Meanwhile, the reconstruction quality of mFPIE, ADMM-FPIE, and AS-FPIE are severely degraded if there is LED shifting, as listed in the Tab. 2 when $d = 0.5$ mm. The performance of these methods will further decrease if the noise signal becomes large for example when $c = 0.25$ and $a = 10^{-2}$. In these cases, our proposed Elfpie still works promising as the LSNR for Elfpie are 32.72 and 35.40, while the LSNR are 1.56, 20.18, and 2.57 for mFPIE, ADMM-FPIE, and AS-FPIE, respectively.

Further increasing the LED shifting range $d$ will lead to the failure of reconstruction of these SOTA methods as some of the LSNRs are less than 0, while the LSNRs of our proposed Elfpie are larger than 30 denoting success reconstruction results.

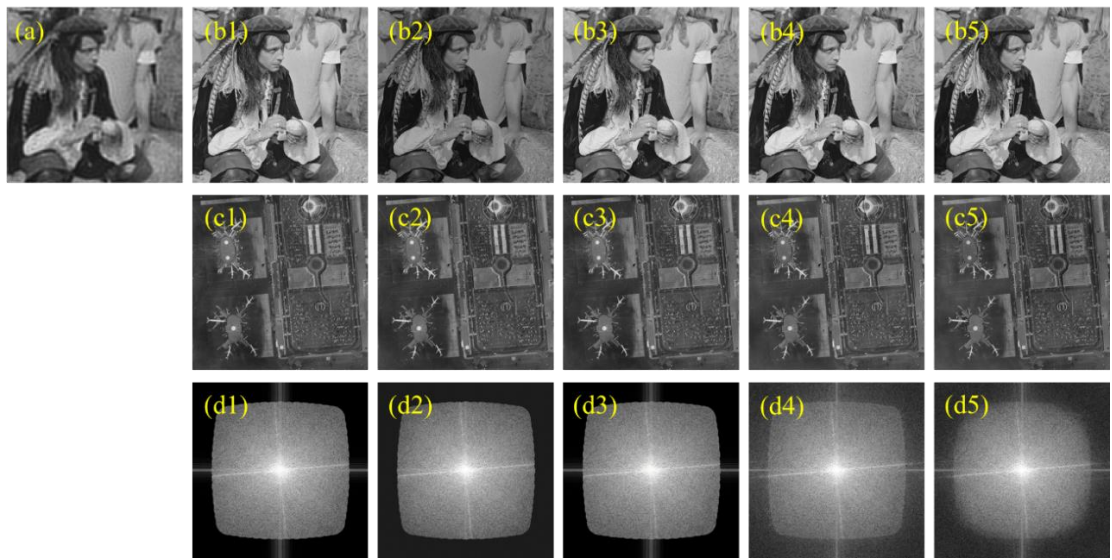

Fig. 4. Reconstruction results for $d = 0$, $c = 0.25$ and $a = 10^{-5}$ Gaussian noise. (a) is the raw image illuminated by the center LED dot. (b1) to (b5) are reconstructed amplitude for mFPIE, ADMM-FPIE, AS-FPIE, and amplitude-based, and intensity-based Elfpie, respectively. (c1) to (c5) are phase component. (d1) to (d5) are Fourier spectra.

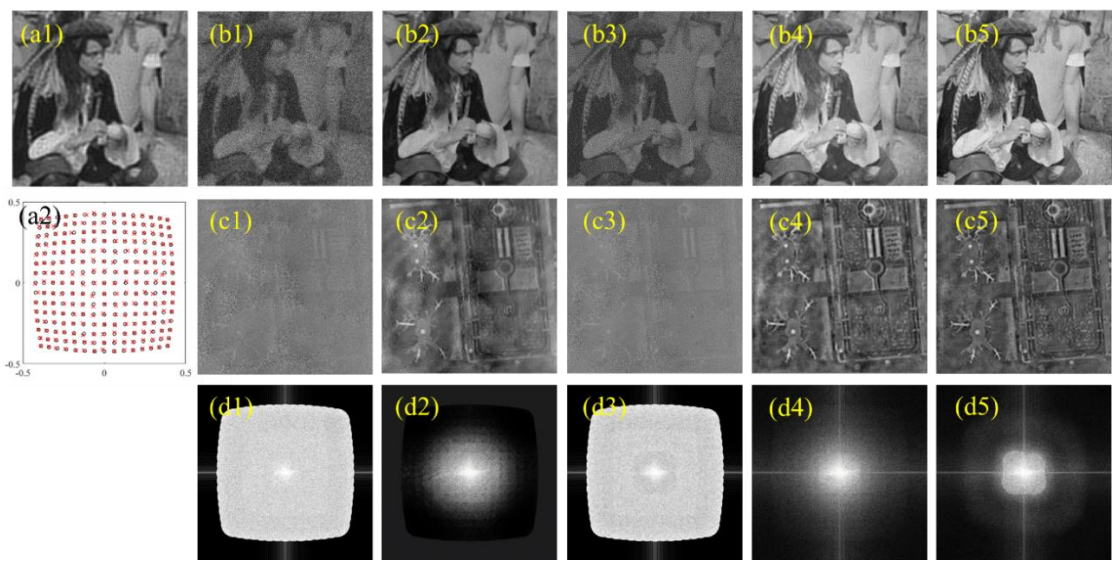

Fig. 5. Reconstruction results for $d = 0.5$, $c = 0.25$ and $a = 10^{-2}$, Gaussian noise. (a1) is the raw image illuminated by the center LED dot. (a2) shows the ideal (red circle) and shifted LED dots (black cross). (b1) to (b5) are reconstructed amplitude for mFPIE, ADMM-FPIE, AS-FPIE, and amplitude-based, and intensity-based Elfpie, respectively. (c1) to (c5) are phase component. (d1) to (d5) are Fourier spectra.

Reconstruction results for $d = 0$, c = 0.25 and $a = 10^{-5}$ are shown in Fig. 4 for visual assessment. According to Figs. 4 (b1) to (b5), all methods recovered promising results for both amplitude and phase components. Interestingly, we found that there are some high-frequency components beyond the recovery region in the Fourier spectrum for Elfpie recovery, as shown in Figs. 4 (d4) and (d5). These high-frequency components are introduced by the Hessian regularization for compensating the noise signals. The Fourier spectrum in Fig. 4 (d5) is not as uniform as that of Figs. 4 (d1) to 4 (d4), this is the property of the intensity-based data fidelity, in which the lower frequency components contain more light energy, and contribute more to the reconstructed spectrum than higher frequency components.

Reconstruction results for $d = 0.5$ mm, c = 0.25 and $a = 10^{-2}$ are shown in Fig. 5 for visual assessment. In such a large noise signal, the mFPIE and AS-FPIE fail to recover the amplitude and phase components. The Elfpie recovers more details for both amplitude [Figs. 5 (b4) and (b5)] and phase [Figs. 5 (c4) and (c5)] than that of ADMM-FPIE [Figs. 5 (b2) and (c2)]. With the Hessian regularization, the Elfpie is able to adaptively recover certain areas of the Fourier spectrum for better reconstruction results. In general, our proposed Elfpie is robust to Gaussian noise 100 times more than mFPIE and AS-FPIE, even if the LED shifting appears.

Apart from Gaussian noise, we have also tested our Elfpie on Poisson noise listed in Tab. 3. In these cases, without the LED shifting, the AS-FPIE obtains the best LSNR score as listed in the first 8 rows when $d = 0$. However, the LSNR of AS-FPIE decreases largely when LED shift is introduced, for example as listed in Tab 3 for $d = 0.5$ mm and 2 mm, while our proposed Elfpie gains better reconstruction results than other methods among different intensity fluctuation, noise strength, and LED shifting levels.

Reconstruction results for $d = 2$ mm, c = 0.75 and Lv1 noise level are shown in Fig. 6 for visual assessment. Due to the mismatch of LED position, the mFPIE, ADMM-FPIE, and AS-FPIE fail to recover the correct amplitude and phase components as shown in Figs. 6 (b1) to 6 (b5), and in Figs. 6 (c1) to 6 (c5). The phase pattern is severely degraded by the unexpected background values as some of the structures (airplanes in the phase pattern) can hardly be observed. On the contrary, the Elfpie still works promising, as both the intensity and phase are recovered and details structure can be observed. For example, the airplanes in the phase pattern can be clearly observed in Figs. 6 (c4) and 6 (c5), although the resolution of the reconstructed phase is not as high as the conditions with no LED position shift.

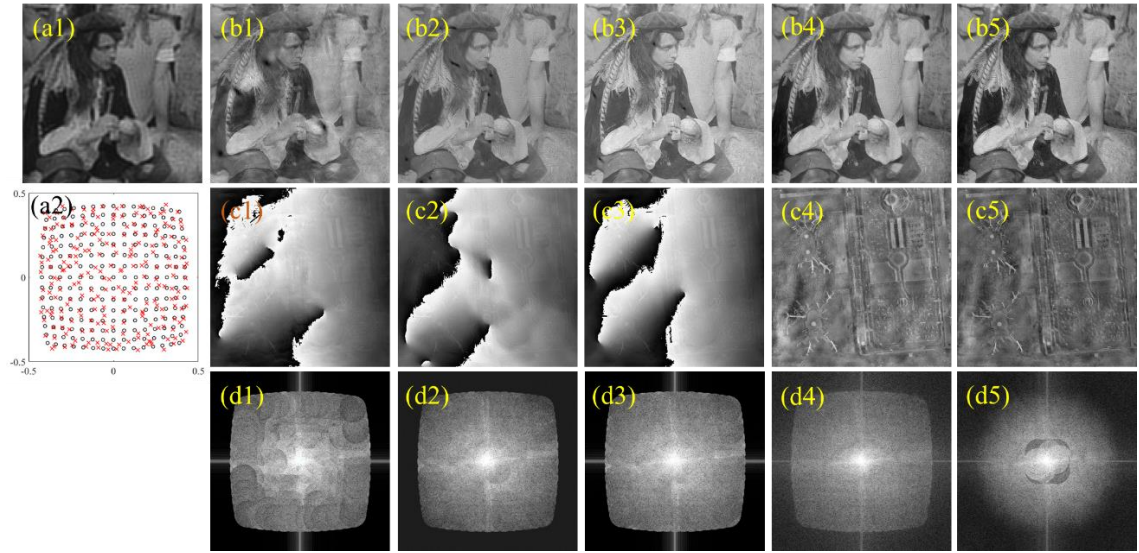

Fig. 6. Reconstruction results for $d = 2$, c = 0.75 and Lv1 noise level. (a1) is the raw image illuminated by the center LED dot. (a2) shows the ideal (red circle) and shifted LED dots (black cross). (b1) to (b5) are reconstructed amplitude for mFPIE, ADMM-FPIE, AS-FPIE, and amplitude-based, and intensity-based Elfpie, respectively. (c1) to (c5) are phase component. (d1) to (d5) are Fourier spectra.

Table 3. Reconstruction quality of different methods for Poisson noise.

| LED position shifting $d$ (mm) | Uneven illumination level $c$ | Noise level | Average LSNR for different methods, **Poisson noise** | | | | |
|---|---|---|---|---|---|---|---|
| | | | mFPM | ADMM | AS | Elfpie_A | Elfpie_I |
| 0 | 0.25 | Lv1 | 52.77 | 55.85 | 57.14 | 50.23 | 48.39 |
| | | Lv2 | 47.21 | 51.22 | 53.26 | 46.49 | 45.05 |
| | | Lv3 | 44.96 | 48.26 | 51.39 | 44.34 | 43.27 |
| | | Lv4 | 41.90 | 44.88 | 48.54 | 41.40 | 41.10 |
| | 0.50 | Lv1 | 44.15 | 49.76 | 50.14 | 47.32 | 46.23 |
| | | Lv2 | 42.18 | 47.77 | 49.35 | 44.82 | 44.03 |
| | | Lv3 | 39.86 | 44.55 | 46.47 | 42.55 | 42.03 |
| | | Lv4 | 37.97 | 42.35 | 44.63 | 40.43 | 40.30 |
| | 0.75 | Lv1 | 34.53 | 42.39 | 41.87 | 43.74 | 43.16 |
| | | Lv2 | 33.02 | 41.18 | 41.22 | 41.71 | 41.47 |
| | | Lv3 | 31.42 | 39.31 | 39.61 | 40.06 | 39.79 |
| | | Lv4 | 31.51 | 37.28 | 38.13 | 38.24 | 37.95 |
| 0.5 | 0.25 | Lv1 | 20.42 | 26.94 | 26.30 | 40.35 | 41.65 |
| | | Lv2 | 18.91 | 28.63 | 25.99 | 40.73 | 41.12 |
| | | Lv3 | 15.52 | 22.50 | 21.35 | 39.14 | 39.10 |
| | | Lv4 | 10.08 | 18.56 | 15.74 | 37.07 | 37.64 |
| | 0.50 | Lv1 | 21.91 | 32.39 | 28.93 | 39.31 | 41.39 |
| | | Lv2 | 14.83 | 21.25 | 21.03 | 39.54 | 39.43 |
| | | Lv3 | 18.99 | 24.30 | 25.11 | 37.97 | 38.22 |
| | | Lv4 | 22.17 | 28.52 | 28.28 | 36.10 | 37.15 |
| | 0.75 | Lv1 | 22.35 | 30.58 | 29.27 | 38.60 | 40.01 |
| | | Lv2 | 20.07 | 29.10 | 27.22 | 38.01 | 38.46 |
| | | Lv3 | 19.82 | 27.95 | 27.23 | 36.80 | 37.20 |
| | | Lv4 | 18.65 | 25.41 | 24.74 | 34.81 | 35.33 |
| 2 | 0.25 | Lv1 | 1.18 | 12.14 | 5.21 | 32.29 | 34.21 |
| | | Lv2 | -0.27 | 13.31 | 4.92 | 32.11 | 33.73 |
| | | Lv3 | 1.06 | 11.18 | 5.01 | 31.05 | 32.22 |
| | | Lv4 | 1.87 | 12.93 | 7.10 | 30.67 | 32.41 |
| | 0.50 | Lv1 | -3.85 | 8.29 | 5.18 | 30.46 | 31.25 |
| | | Lv2 | 1.03 | 15.95 | 5.98 | 31.34 | 32.93 |
| | | Lv3 | 2.66 | 11.87 | 8.04 | 31.26 | 32.92 |
| | | Lv4 | 1.21 | 12.23 | 8.41 | 29.80 | 31.01 |
| | 0.75 | Lv1 | 0.23 | 11.56 | 5.89 | 31.38 | 33.26 |
| | | Lv2 | 4.83 | 14.14 | 10.60 | 31.38 | 32.69 |
| | | Lv3 | 5.38 | 16.05 | 15.36 | 30.66 | 32.25 |
| | | Lv4 | -1.16 | 9.09 | 4.72 | 29.39 | 30.33 |

Table 4. Reconstruction quality of different methods for SNP noise.

| LED position shifting $d$ (mm) | Uneven illumination level $c$ | Noise level $a$ | Average LSNR for different methods, **SNP noise** | | | | |
|---|---|---|---|---|---|---|---|
| | | | mFPM | ADMM | AS | Elfpie_A | Elfpie_I |
| 0 | 0.25 | 0.001 | 29.65 | 51.27 | 38.99 | 51.38 | 43.17 |
| | | 0.01 | 10.55 | 42.76 | 21.47 | 41.63 | 39.16 |
| | | 0.1 | -1.85 | 33.41 | 5.98 | 43.16 | 33.09 |
| | | 0.2 | -6.03 | 29.91 | -0.77 | 34.60 | 30.27 |
| | 0.50 | 0.001 | 27.61 | 45.11 | 36.53 | 46.36 | 40.46 |
| | | 0.01 | 9.86 | 39.42 | 20.60 | 39.08 | 37.52 |
| | | 0.1 | -2.12 | 32.25 | 5.51 | 41.96 | 32.28 |
| | | 0.2 | -6.16 | 29.25 | -1.09 | 34.13 | 29.98 |
| | 0.75 | 0.001 | 24.30 | 38.30 | 32.57 | 42.93 | 38.11 |
| | | 0.01 | 8.96 | 35.28 | 19.52 | 36.56 | 35.54 |
| | | 0.1 | -2.65 | 30.03 | 4.72 | 39.61 | 31.24 |
| | | 0.2 | -6.37 | 27.54 | -1.48 | 33.48 | 29.04 |
| 0.5 | 0.25 | 0.001 | 14.83 | 27.52 | 20.41 | 40.64 | 38.48 |
| | | 0.01 | 7.61 | 32.32 | 16.64 | 37.26 | 35.97 |
| | | 0.1 | -4.83 | 20.63 | 0.84 | 33.79 | 30.88 |
| | | 0.2 | -7.50 | 19.90 | -3.10 | 30.64 | 30.00 |
| | 0.50 | 0.001 | 15.31 | 27.06 | 21.49 | 38.30 | 37.10 |
| | | 0.01 | 5.84 | 26.63 | 14.27 | 35.79 | 34.52 |
| | | 0.1 | -5.25 | 19.43 | 0.14 | 34.23 | 30.81 |
| | | 0.2 | -8.57 | 12.61 | -4.85 | 30.67 | 28.68 |
| | 0.75 | 0.001 | 13.09 | 21.77 | 18.27 | 36.32 | 34.86 |
| | | 0.01 | 2.23 | 21.68 | 9.78 | 34.02 | 32.83 |
| | | 0.1 | -6.15 | 14.71 | -1.47 | 33.82 | 30.14 |
| | | 0.2 | -8.04 | 18.75 | -3.83 | 30.78 | 28.55 |
| 2 | 0.25 | 0.001 | -0.70 | 14.13 | 5.92 | 31.46 | 32.64 |
| | | 0.01 | -4.19 | 13.19 | 2.65 | 30.91 | 31.25 |
| | | 0.1 | -7.64 | 11.06 | -3.29 | 30.09 | 29.45 |
| | | 0.2 | -9.89 | 8.29 | -6.79 | 29.16 | 28.54 |
| | 0.50 | 0.001 | -1.52 | 8.48 | 2.98 | 30.29 | 30.25 |
| | | 0.01 | -0.76 | 15.10 | 4.70 | 30.64 | 30.43 |
| | | 0.1 | -6.68 | 11.58 | -2.38 | 30.11 | 29.25 |
| | | 0.2 | -9.55 | 9.42 | -6.33 | 29.32 | 28.46 |
| | 0.75 | 0.001 | 1.95 | 10.41 | 7.04 | 30.37 | 30.18 |
| | | 0.01 | -4.01 | 8.44 | 0.47 | 30.54 | 29.71 |
| | | 0.1 | -7.39 | 7.24 | -3.71 | 29.02 | 27.95 |
| | | 0.2 | -9.50 | 6.07 | -6.48 | 28.16 | 27.27 |

The Elfpie is also robust to SNP noises. As listed in Tab. 4. The Elfpie is robust to very large SNP noise, for example, $a = 0.2$ where the noise strength is 95571.99 %. This property is attributed to the combination of $L_1$-norm in the data fidelity and the Hessian regularization. When a = 0.1, the LSNRs for both mFPIE and AS-FPIE are less than 0 denoting the failure of reconstruction. Results for $d = 0$, c = 0.25 and $a = 0.2$ are shown in Fig. 7 for visual assessment. Accordingly, the Elfpie succeeds to reconstruct the high-quality amplitude and phase component, while all other methods fail.

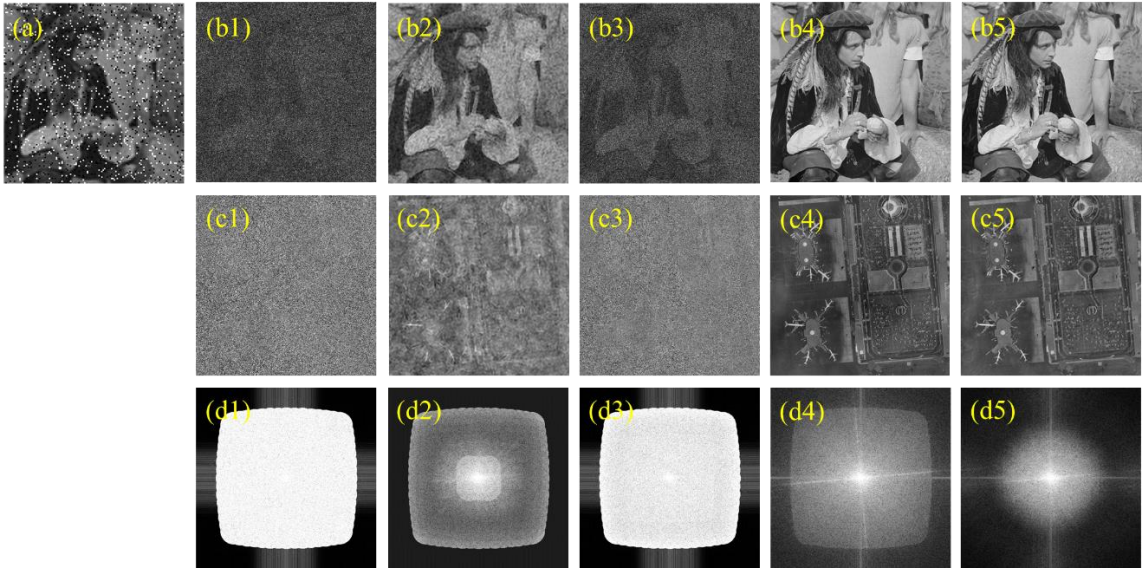

Fig. 7. Reconstruction results for $d = 0$, $c = 0.25$ and $a = 0.2$ for SNP noise. (a) is the raw image illuminated by the center LED dot. (b1) to (b5) are reconstructed amplitude for mFPIE, ADMM-FPIE, AS-FPIE, and amplitude-based, and intensity-based Elfpie, respectively. (c1) to (c5) are phase component. (d1) to (d5) are Fourier spectra.

### 4.3 Bypassing the vignetting effect

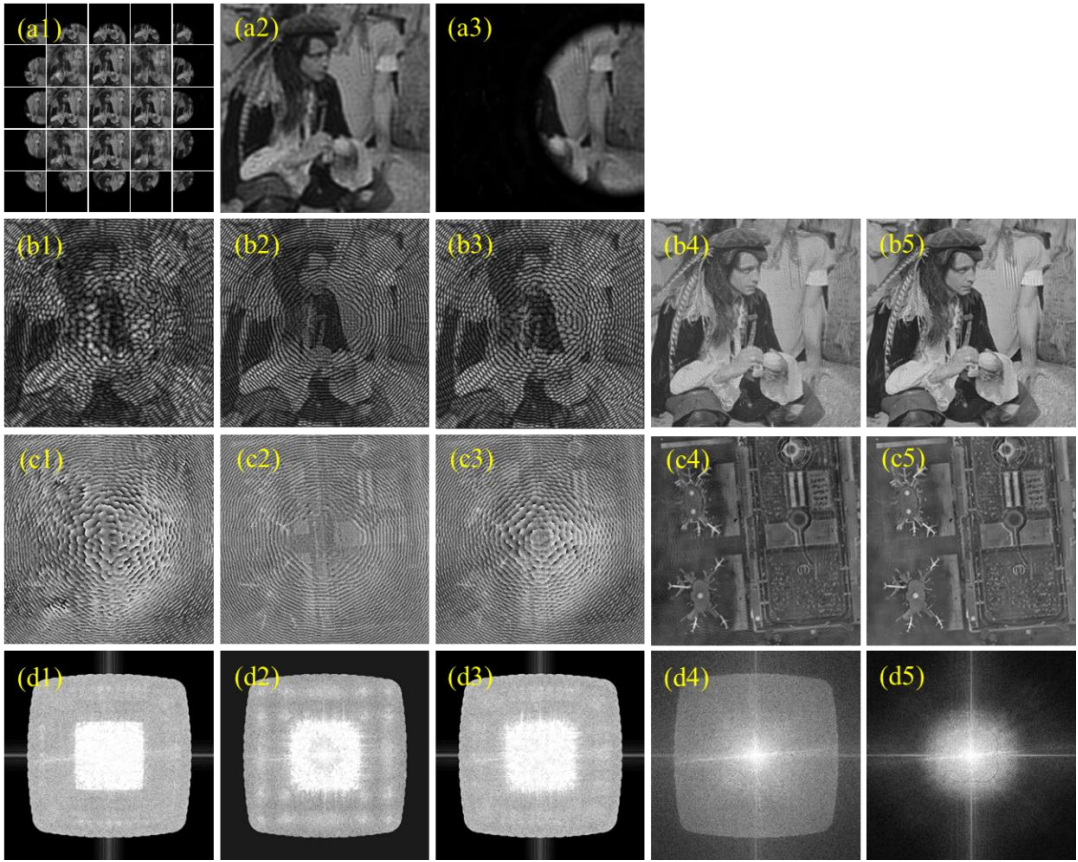

Fig. 8, experiment on FPM images with vignetting effect. (a1) first 25 images from the FPM data cube. (a2) bright-field image. (a3) dark-field image with vignetting effect. (b1) to (b5) are reconstruction amplitude using mFPIE, ADMM-FPIE, AS-FPIE, and Elfpie with amplitude and intensity data fidelity, respectively. (c1) to (c5) phase pattern. (d1) to (d5) are Fourier spectrums for (b1) to (b5).

Our proposed Elfpie can also bypass the vignetting effect without any modification to the raw data. To show this, we generate the simulated data with vignetting effect as shown in Fig. 8. Figure 8 (a1) shows the first 25 images in the FPM data cube. The Gaussian noise strength is $a = 0.0001$. Not LED shifting is added. and Figs. 8 (a2) and 8 (a3) show the central image, and one of the vignetting images where the image is partially in a bright field and partially in a dark field.

Figure 8 (b1) to (b5) shows the reconstruction results direct from the problematic data. Accordingly, all the SOTA methods fail to reconstruct the high-quality data as the amplitude and phase patterns are severely degraded by the interfering artifacts. Since vignetting images contain unexpected low-frequency information of high-energy, direct recovering using vignetting images introduces the high-energy low-frequency information to the Fourier spectrum as shown in Figs. 8 (d1) to (d3), where a large area near the center of Fourier spectrum is bright, and further degraded the reconstruction quality.

While in our proposed Elfpie, we use the spatial gradient to recover the image. By taking the image gradient to the vignetting images the low-frequency component is removed, and the Fourier spectrum can be recovered more uniformly, as shown in Figs. 8 (d4) and 8 (d5). In such a manner, we are able to bypass the vignetting effect, and obtain promising amplitude, as shown in Figs. 8 (b4) and 8 (b5), and phase as shown in Figs. 8 (c4) and 8 (c5) reconstruction results

According to our simulation tests, our proposed method Elfpie is robust to Gaussian noise, Poisson noise, and SNP noise from weak to strong levels. The reconstruction quality of Elfpie is also not sensitive to the LED position shifting up to 2 mm. In addition, we found that the intensity-based data fidelity in our Elfpie works better than the amplitude-based Elfpie in more cases such as Poisson and Gaussian noise.

## 5. Experimental study

### 5.1 *Testing on public datasets*

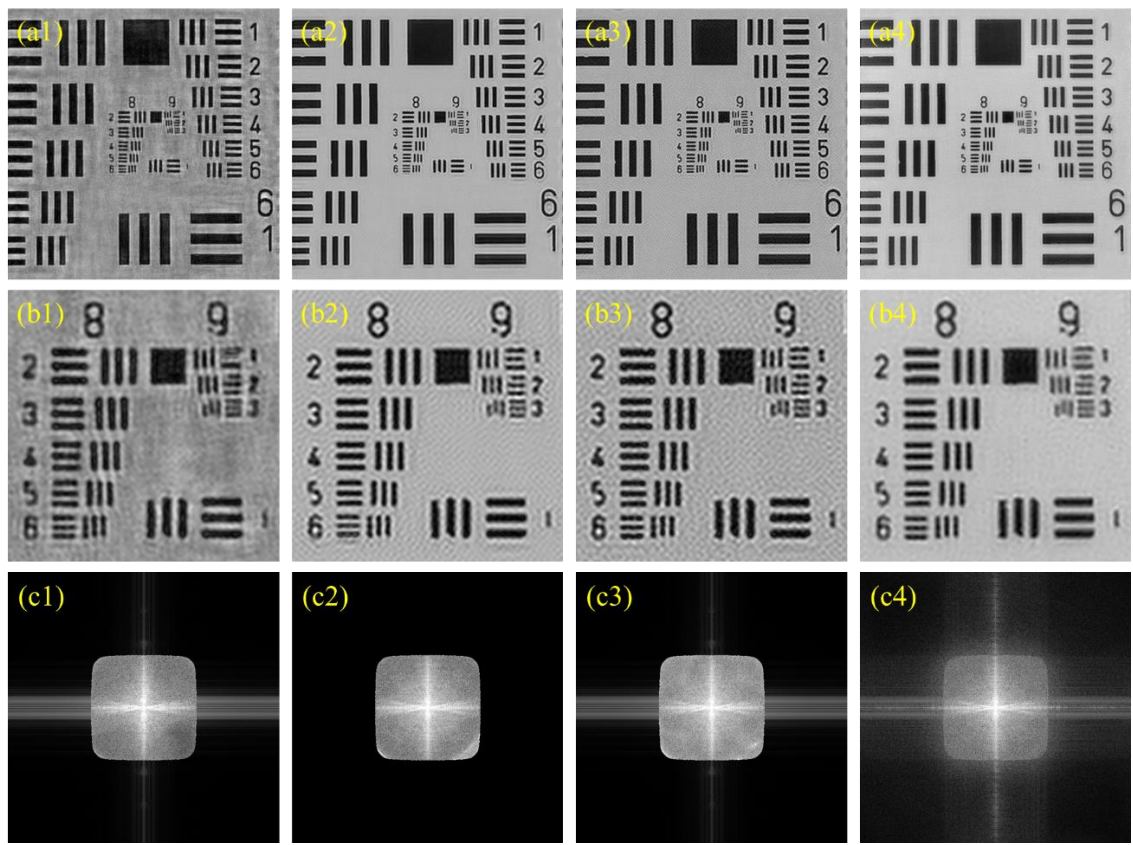

Fig. 9, Reconstruction results for Zuo's dataset. (a1) to (a4) are amplitude patterns for mFPIE, ADMM-FPIE, AS-FPIE, and Elfpie. (b1) to (b4) are zoomed-in parts for groups 8 and 9. (c1) to (c4) are Fourier spectra.

We perform the experimental study on the public FPM datasets, and compare our Elfpie against SOTA methods. First, we test the Elfpie in a USAF target from Zuo's dataset which can be downloaded from

([https://www.scilaboratory.com/code.html](https://www.scilaboratory.com/code.html)). As shown in Figs. 9 (a1) to 9 (a3) and Figs. 9 (b1) to 9 (b3), the results of mFPIE, ADMM-FPIE, and AS-FPIE are degraded by the noise signals and these images have obvious ringing effect which is caused by the noise pixels in the raw data cube. The noise effect in ADMM-FPIE shown in Figs. 9 (a2) and 9 (b2) is less than that in Figs. 9 (a3) and 9 (b3) for AS-FPIE. While our proposed Elfpie has better-denoised results as the background, shown in Figs. 9 (a4) and 9 (b4), is more flattened than that of other methods.

We further test the Elfpie in Zheng's dataset, which can be downloaded from [https://github.com/SmartImagingLabUConn/Fourier-Ptychography](https://github.com/SmartImagingLabUConn/Fourier-Ptychography), on both blood smear in Fig. 10 and MouseKidney slice in Fig. 11. As shown in Figs. 10 (b1) to 10 (b4), our Elfpie has no boundary effect in the recovered phase pattern while others have. As shown in Fig. 11 (a4) and 11 (b4), both the amplitude and phase pattern of the Elfpie have better recovery quality than others in Figs. 11 (a1) to 11 (b3).

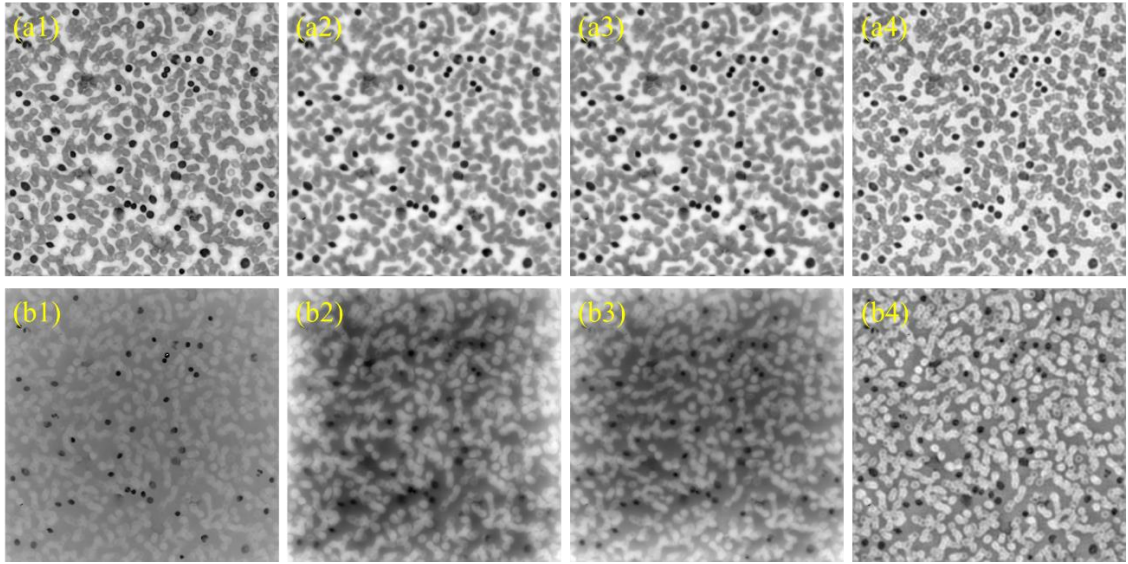

Fig. 10, Reconstruction results on 'bloodsmear_green' data. (a1) to (a4) are amplitude for mFPIE, ADMM-FPIE, AS-FPIE, and Elfpie. (b1) to (b4) are phase patterns.

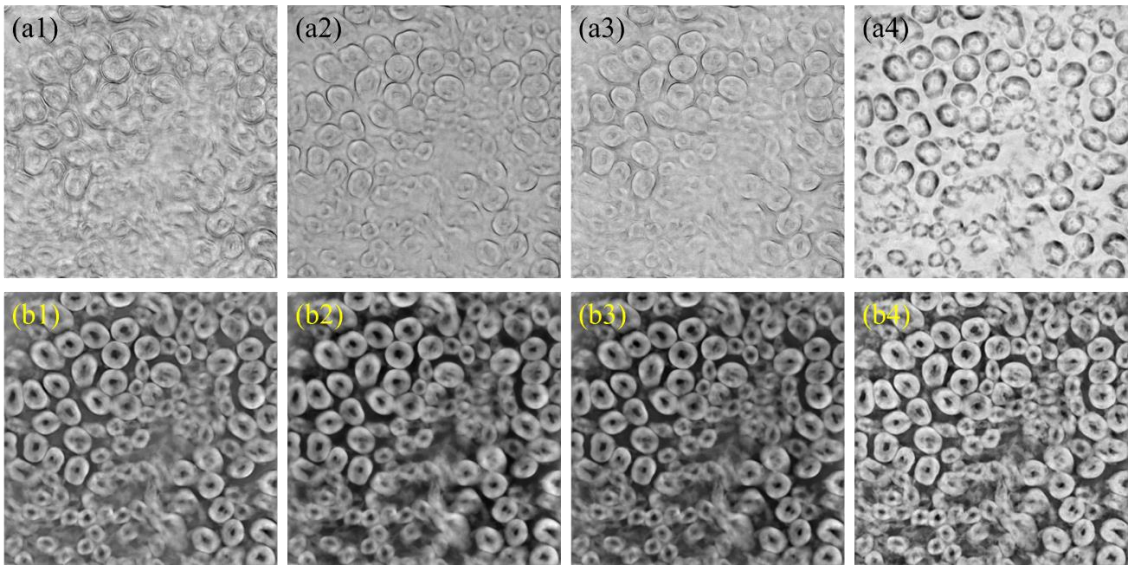

Fig. 11, Reconstruction results on 'MonkeyKidney_green' data. (a1) to (a4) are amplitude for mFPIE, ADMM-FPIE, AS-FPIE, and Elfpie. (b1) to (b4) are phase patterns.

*5.2 Testing on the private dataset*

We also test the Elfpie on our private dataset which is collected by a simply constructed FPM platform. The camera pixel size is 3.45 μm, and the objective lens is NA = 0.1, ×3. The LED is placed 180 mm away from the sample, and the distance between two adjacent LED dots is 8 mm. The illumination wavelength is 522 nm. A total of 64 LEDs is used. The raw data is of 260 by 260 pixels as shown in Fig. 12 (a1). Due to the simple structure of the FPM platform, the raw image has the vignetting effect and severe Gaussian noise as shown in Fig. 12 (a2). Potential LED position misaligning also exists. Pupil correction is turned off. 50 iterations for all methods.

The reconstructed amplitude patterns are shown in Figs. 12 (b1) to 12 (b4). Figures 12 (c1) to 12 (c4) show the zoomed-in part for group 8 to group 9. Accordingly, the mFPIE, ADMM-FPIE, and AS-FPIE fail to recover the pattern as their results are corrupted by heavy artefacts. The stripe pattern is completely merged into the problematic pixels. While the Elfpie can recover the high-resolution, noise-free images as shown in Figs. 12 (b4) and 12 (c4). The stripe pattern for Group 9 element 1 can be clearly observed.

Moreover, as shown in Figs. 12 (d1) to (d4), the Fourier spectrums for these three methods are severely non-uniformed, as rapid changes in the power spectrum can be observed for example in Figs. 12 (d1) and 12 (d2). While for Elfpie, as shown in Fig. 12(d4), the Fourier spectrum is more uniform and contains more information which suppresses the noise signals during the reconstruction.

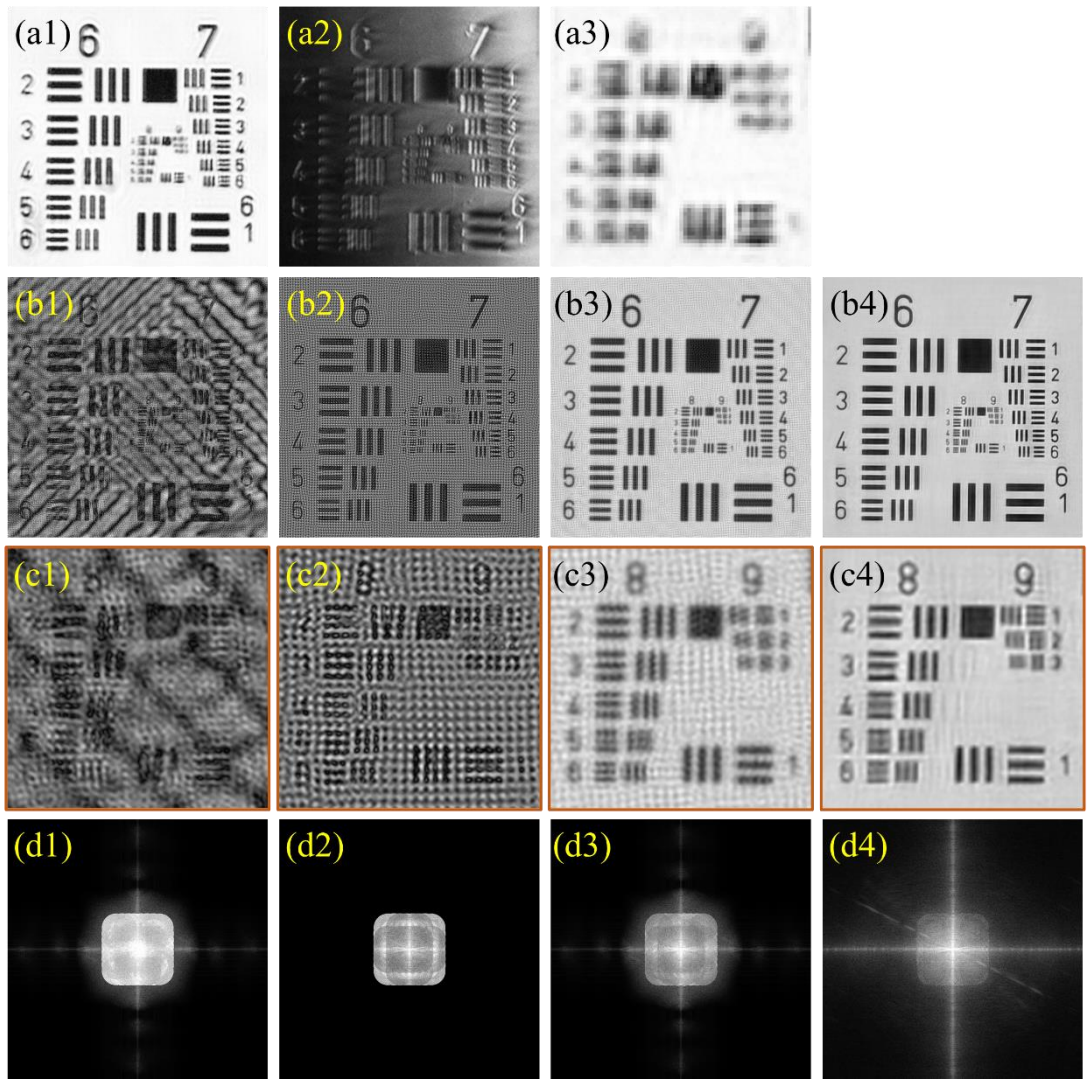

Fig. 12, reconstruction results on USAF target. (a1) bright-field raw image. (a2) dark-field raw image. (a3) zoomed-in image for groups 8 to 9 in (a1). (b1) to (b4) are reconstruction amplitudes using mFPIE, ADMM-FPIE, AS-FPIE, and Elfpie, respectively. (c1) to (c4) are zoomed-in images for groups 8 to 9. (d1) to (d4) are Fourier spectrums for (b1) to (b4).

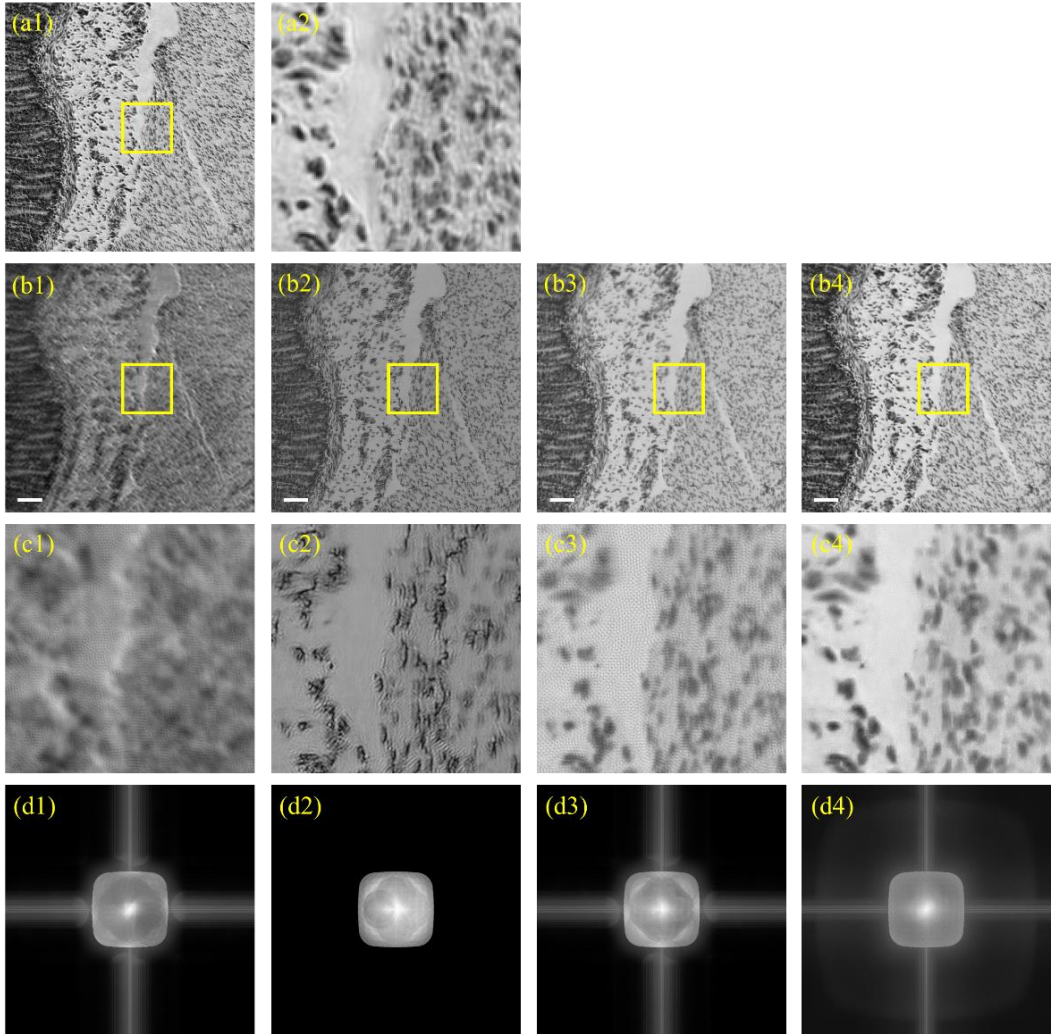

Fig. 13, reconstruction results on intestinal epithelial slice. (a1) bright-field raw image. (a2) zoomed-in image for the yellow box in (a1). (b1) to (b4) are reconstruction amplitudes using mFPIE, ADMM-FPIE, AS-FPIE, and Elfpie, respectively. (c1) to (c4) are zoomed-in images for area in yellow boxes. (d1) to (d4) are Fourier spectrums for (b1) to (b4). Scale bar is 50 μm.

We perform another group of experiment on intestinal epithelial slice. The camera pixel size is 6.65 μm, and the objective lens is NA = 0.25, ×10. The LED is placed 80 mm away from the sample. The illumination wavelength is 630 nm. A total of 64 LEDs is used. The raw data is of 768 by 768 pixels as shown in Fig. 13 (a1). 30 iterations for all methods. Again, the mFPIE fails to recover the sample amplitude due to the low signal-to-noise ratio of the raw image. According to the results in Figs. 13 (c1) to 13 (c4), our Elfpie obtains high quality reconstruction result, in which the noise signal is efficiently suppressed, and the small structures can be correctly recovered.

## 6. Discussion

### 6.1 Elfpie cost function

In this section, we analyze the function of each term in the cos function of Elfpie. First, our Elfpie uses the image spatial gradient, sharing the similar idea of the variation Retinex theory [42-45], for FPM reconstruction. By doing this, it is able to bypass the impact of background intensity, including vignetting effect. This is because by taking the gradient, we efficiently remove the low-frequency component of the raw images. Although the gradient operators may also increase the noise signal, the Hessian penalty is strong enough to suppress almost all kinds of noises including Gaussian, Poisson [46]. Together with the $L_1$-norm fidelity, the Elfpie is also robust to the SNP noises [47, 48] while other FPIEs fail. Furthermore, by using the edge information for the recovery, the Elfpie has stronger tolerance to the LED position shifting than previous FPIEs.

The Hessian regularization is applied to both the amplitude and phase of the high-resolution image. The Hessian regularization uses the image's second-order spatial gradient and won't cause a staircase-like effect as the TV regularization does. Unlike other FPIE with TV regularization where the TV denoising is not embedded into the gradient descent process [30, 31], our Hessian regularization is directly incorporated into the cos function, and the closed-form gradient is provided. In such a manner, our Hessian denoising yields stable convergent properties.

*6.2 Implementation efficiency*

The Elfpie has 5 hyper-parameters that need to be manually adjusted: α and β for Hessian regularization, and r1, r2, and initial step size for the optimizer. Fortunately, tuning these hyperparameters is not time-consuming and is way easier than tuning the layout of optical system as the execution speed of Elfpie is fast. Moreover, r1¬ = 0.9, and r2 = 0.999 are testing in our massive simulations for good reconstruction quality, and α and β are given by Eq. (19), which can be regarded as references value for fine-tuning.

Since our Elfpie can be performed in parallel, the speed of our Elfpie is significantly faster than other SOTA FPIEs. The following Tab. 5 shows the average execution time for mFPIE, ADMM-FPIE, AS-FPIE, and Elfpie when parallel computation is turned on (4 Threads). The testing image is 141 by 141 pixels, a total of 225 images in FPM the data cube. The code is running on MATLAB 2022, with Intel(R) Core(TM) i5-8265U CPU @ 1.60GHz, 8.00 GB RAM.

**Table 5. Reconstruction speed for different methods.**

| Methods | mFPIE (50 iterations) | ADMM-FPIE (50 iterations) | AS-FPIE (30 iterations) | Elfpie (50 iterations) |
|---|---|---|---|---|
| Time (s) | 100.35 | 76.09 | 34.47 | 23.4 |

Our proposed optimizer share the advatanges of AdaBelief and AdaDelta, it converges fast and can skip the local minima. We compare the proposed optimizer against other off-the-shelf including the SGD (stochastic gradient descent), the NAdam, and the AdaBelief based on the experiment in Fig. 12. The change of fidelity loss with respect to iterations are shown in Fig. 13 (a). Accordingly, the SGD converges the fastest but falls into local minima as shown in the black curve in Fig. 13 (a) after 60 iterations. While our proposed optimizer converges faster than that of the NAdam [49], and AdaBelief, with the same parameters, as shown in the red curve in Fig. 13 (a). By increasing the step-size, the proposed optimizer converges faster, but to the same minimal as shown in the magenta curve as shown in Fig. 13 (b). Note that the fidelity loss enters lower energe than the black curve.

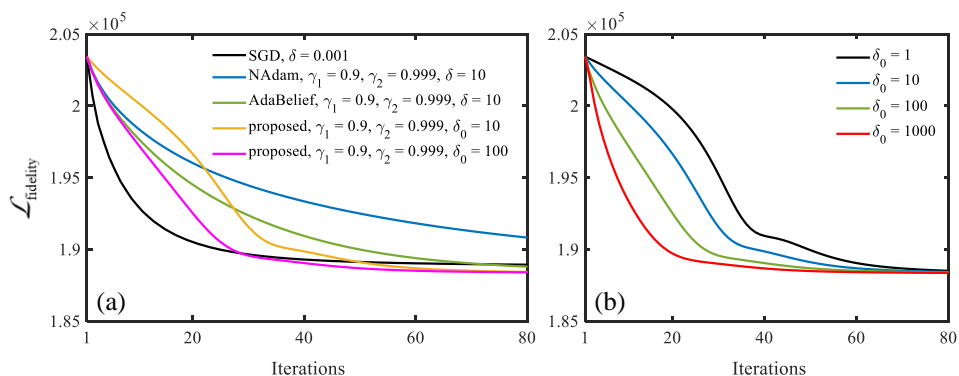

Fig. 13, Elfpie convergence properties. (a) convergence behavior of different optimizers. (b) convergence behavior of different initial learning rate.

*6.3 Amplitude-based fidelity as a special case of Gamma correction*

Although the intensity-based fidelity works better than amplitude-based fidelity according to the simulation results, we found that this conclusion is not always true as sometimes the amplitude-based fidelity works better than that of intensity-based fidelity in practical experiments. As present, we are unable to confirm which one is the best since both are mathematically correct, but we are able to consider the function $g$ in Eq. (3) to be the gamma correction function where $g(x) = x^{\gamma}$. As such, the amplitude-based fidelity can be regarded as a special case of gamma correction where $\gamma = 0.5$.

We have also tested the reconstruction results with different value of $\gamma$, as shown in Fig. 14 (a1) to 14 (a4). The Hessian regularization is turned off. Accordingly, we found that, in this experiment, small value of $\gamma$ is able to suppress the noise signals, while the reconstructed image is over smoothed as shown in Fig. 14 (a1) for $\gamma = 0.01$. The high-frequency component is less recovered which can be seen from Eq. 14 (c1), where the brightness of high-frequency is less than that in Fig. 14 (c2) or 14 (c3).

When $\gamma = 0.875$, as shown in Fig. 14 (a4), the reconstructed image is corrupted by noise signals, but the sharpness is increased. We have further tried $g(x) = \log(x+1)$ as shown in Fig. 14 (a5), in this case, the sharpness of the reconstructed image is maintained, and noise signals is also suppressed compared to Fig. 14 (a4). But the convergence is not good since the stripes in group 9 element 1 are not well reconstructed.

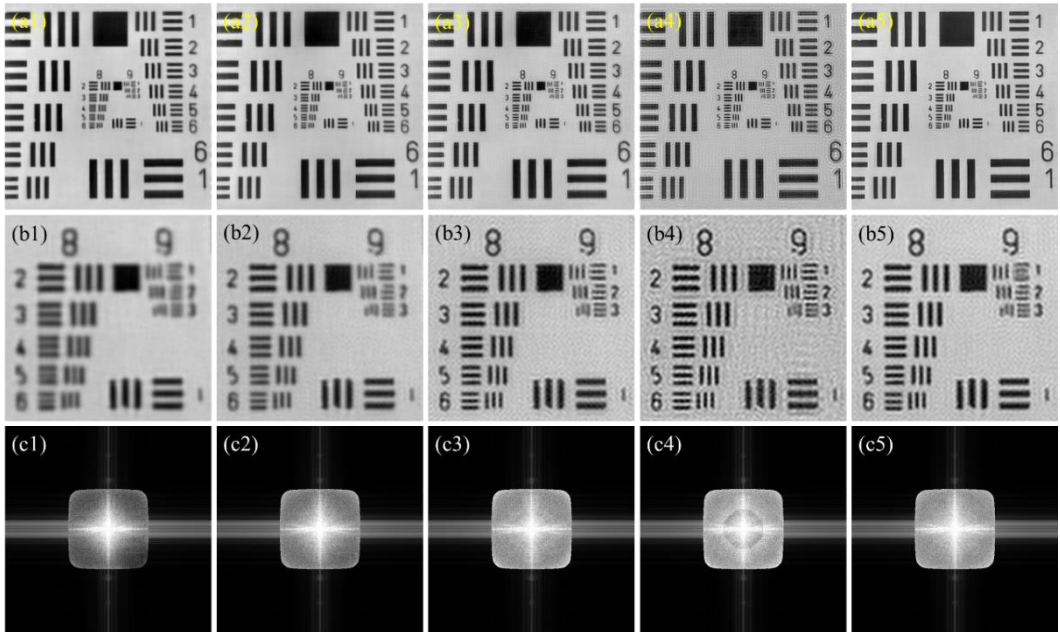

Fig. 14, Elfpie reconstruction using different scaling function $g$. Hessian regularization is turned off, $\alpha = \beta = 0$.
(a1) $g(x) = x^{0.01}$; (a2) $g(x) = x^{0.125}$; (a3) $g(x) = x^{0.375}$; (a4) $g(x) = x^{0.875}$; (a5) $g(x) = \ln(x+1)$. (b1) to (b5) are zoomed-in for groups 8 and 9. (c1) to (c5) are Fourier spectra. 30 iterations.

## 7. Concluding remarks

In this research, we propose the Elfpie for error-laxity Fourier ptychographic reconstruction. The Elfpie uses image's first-order gradient to form the data fidelity term, which is robust to LED misaligning, and problematic illuminations. The second-order total variation regularization (Hessian regularization) is also incorporated in the Elfpie to suppress the noise signal including Gaussian, Poisson, and SNP noise. The Elfpie achieves high-quality FPM reconstruction and is more robust to system's error than other SOTA reconstruction algorithm. It largely decreases the system calibration requirement. As such, the FPM technique can become more user friendly and robust to different working environments. Future works including extending Elfpie to 3D FPM, lensless FPM, and other related applications. Investigation of different penalty function or data fidelity function (including $\Theta$, $g$, and distance function) can also potentially optimize the Elfpie.

**Declaration of Competing Interest**
The authors declare that they have no known competing financial interests or personal relationships that could have appeared to influence the work reported in this paper.

21 / 23